\documentclass{amsart}

\usepackage{amsmath,amsthm,amssymb,amsfonts,mathrsfs}
\usepackage{enumerate}
\usepackage{geometry}
\usepackage[colorlinks]{hyperref}
\usepackage[nameinlink]{cleveref}

\usepackage{url}
\usepackage{ifthen}
\usepackage[all]{xy} \xyoption{2cell} \UseAllTwocells
\usepackage[mathscr]{euscript}
\usepackage[T1]{fontenc}
\usepackage[utf8]{inputenc}   % For UTF-8 input encoding
\usepackage{refcount}
\usepackage{zref-user}
\usepackage{zref-abspage}
\usepackage{tikz,tikz-cd}
\usepackage{comment}
\numberwithin{equation}{section}

\theoremstyle{plain}
\newtheorem{thm_}[equation]{Theorem}
\newtheorem{lemma_}[equation]{Lemma}
\newtheorem{prop_}[equation]{Proposition}
\newtheorem{cor_}[equation]{Corollary}
\newtheorem{eg_}[equation]{Example}

\newtheorem{claim_}[equation]{Claim}

\theoremstyle{definition}
\newtheorem{thmu_}[equation]{Theorem}
\newtheorem*{thmus_}{Theorem}
\newtheorem{propu_}[equation]{Proposition}
\newtheorem*{propus_}{Proposition}
\newtheorem{coru_}[equation]{Corollary}
\newtheorem*{corus_}{Corollary}
\newtheorem{lemu_}[equation]{Lemma}
\newtheorem*{lemus_}{Lemma}
\newtheorem{egu_}[equation]{Example}
\newtheorem*{egus_}{Example}
\newtheorem{def_}[equation]{Definition}
\newtheorem*{defs_}{Definition}
\newtheorem{rk_}[equation]{Remark}
\newtheorem*{rks_}{Remark}
\newtheorem{ex_}[equation]{Remark}
\newtheorem*{exs_}{Remark}
\newtheorem{constr_}[equation]{Construction}
\newtheorem*{constrs_}{Construction}
\newtheorem{nota_}[equation]{Notation}
\newtheorem*{notas_}{Notation}

\newcommand{\thm}[1]{\begin{thm_}#1\end{thm_}}
\newcommand{\thmu}[1]{\begin{thmu_}#1\end{thmu_}}

\newcommand{\lemm}[1]{\begin{lemma_}#1\end{lemma_}}

\newcommand{\prop}[1]{\begin{prop_}#1\end{prop_}}

\newcommand{\defi}[1]{\begin{def_}#1\end{def_}}

\newcommand{\rk}[1]{\begin{rk_}#1\end{rk_}}

\newcommand{\cor}[1]{\begin{cor_}#1\end{cor_}}

\newcommand{\pf}[1]{\begin{proof}#1\end{proof}}

\DeclareMathOperator{\GL}{GL}
\DeclareMathOperator{\PGL}{PGL}
\DeclareMathOperator{\SL}{SL}

\DeclareMathOperator{\Gal}{Gal}
\DeclareMathOperator{\Hom}{Hom}

\DeclareMathOperator{\Spec}{Spec}

\DeclareMathOperator{\im}{im}
\DeclareMathOperator{\diag}{diag}
\DeclareMathOperator{\Ob}{Ob}

\DeclareMathOperator{\Br}{Br}

\DeclareMathOperator{\Pic}{Pic}

\DeclareMathOperator{\inv}{inv}

\DeclareMathOperator{\Map}{Map}
\newcommand{\Set}{{\sS\tu{et}}}

\newcommand{\Sch}{{\sS\tu{ch}}}
\newcommand{\Esp}{{\sE\tu{sp}}}
\newcommand{\Chp}{{\sC\tu{hp}}}

\newcommand{\CC}{\mathbb C}

\newcommand{\QQ}{\mathbb Q}

\newcommand{\ZZ}{\mathbb Z}
\newcommand{\bfA}{\mathbf A}%
\newcommand{\bfG}{\mathbf G}%

\newcommand{\sC}{\mathscr C}
\newcommand{\sD}{\mathscr D}
\newcommand{\sE}{\mathscr E}

\newcommand{\sS}{\mathscr S}

\newcommand{\cO}{\mathcal O}%

\newcommand{\cV}{\mathcal V}

\newcommand{\cX}{\mathcal X}
\newcommand{\cY}{\mathcal Y}
\newcommand{\cZ}{\mathcal Z}%

\newcommand{\s}{\sigma}

\newcommand{\tm}{\times}%

\newcommand{\ol}{\overline}

\newcommand{\ra}{\rightarrow}

\newcommand{\xra}{\xrightarrow}

\newcommand{\mpt}{\mapsto}

\newcommand{\is}[2]{\xymatrix@-4mm{#1 \ar[r]^-{\sim} & #2 }}
\newcommand{\mis}[2]{\xymatrix@-2mm{#1 \ar[r]^-{\sim} & #2 }}

\newcommand{\dra}[4]{\xymatrix@-4mm{#1 \ar@<.5ex>[r]^-{#3} \ar@<-.5ex>[r]_-{#4}& #2 }}  % double ra
\newcommand{\era}[5]{\xymatrix@-4mm{#1 \ar[r] &#2 \ar@<.5ex>[r]^-{#4} \ar@<-.5ex>[r]_-{#5}& #3 }}  % equalizer
\newcommand{\texop}[2]{{\texorpdfstring{#1}{#2}}}
\newcommand{\tc}[1]{{\text{\textcircled{#1}}}}

\newcommand{\tu}[1]{\text{\upshape #1}}

\newcommand{\op}{{op}}
\newcommand{\et}{\tu{\'et}}

\newcommand{\fppf}{\tu{fppf}}

\DeclareFontFamily{U}{wncy}{}
\DeclareFontShape{U}{wncy}{m}{n}{%
   <5>wncyr5%
   <6>wncyr6%
   <7>wncyr7%
   <8>wncyr8%
   <9>wncyr9%
   <10>wncyr10%
   <11>wncyr10%
   <12>wncyr6%
   <14>wncyr7%
   <17>wncyr8%
   <20>wncyr10%
   <25>wncyr10}{}
\DeclareMathAlphabet{\cyrille}{U}{wncy}{m}{n}
\def\Sha{\cyrille X}

\newcommand{\Tors}{\mathsf{Tors}}

\newcommand{\eq}[1]{\begin{equation}#1\end{equation}}
\newcommand{\eqn}[1]{\begin{equation*}#1\end{equation*}}
\newcommand{\ga}[1]{\begin{gather}#1\end{gather}}

\newcommand{\aln}[1]{\begin{align*}#1\end{align*}}

\newcommand{\enmt}[1]{\begin{enumerate}#1\end{enumerate}}

\newcommand{\aci}[1]{\ar@{^(->}[#1]|-{/}}
\newcommand{\coaci}[1]{\ar@{_(->}[#1]|-{/}}
\newcommand{\aoi}[1]{\ar@{^(->}[#1]|-{\circ}}
\newcommand{\coaoi}[1]{\ar@{_(->}[#1]|-{\circ}}
\makeatletter
\def\citet@url@sp{https://stacks.math.columbia.edu/}
\def\citet@bib@sp{stacks-project}
\def\citet@url@kd{https://kerodon.net/}
\def\citet@bib@kd{kerodon}
\newcommand{\citet@tag}[2]{\href{#2tag/#1}{#1}}
\newcommand{\citet@taglist}[2]{%
 \def\@citet@e{}%
 \def\@citet@tag@n{0}% count first
 \@for\@citet@tag:=#1\do{%
  \edef\@citet@tag@n{\the\numexpr\@citet@tag@n + 1}%
 }%
 \def\@citet@tags{%
  \def\@citet@tag@i{0}%
  \@for\@citet@tag:=#1\do{%
   \edef\@citet@tag@i{\the\numexpr\@citet@tag@i + 1}%
   \ifthenelse{\@citet@tag@i > 1}{% multiple tags
    \ifthenelse{\@citet@tag@i = \@citet@tag@n}{% the last loop
     \citet@seplast%
    }{%
     \citet@sep%
    }%
   }{}%
   \citet@entry{\@citet@tag}{#2}% fixme: how to allow spaces around the input tags
  }%
 }%
 \ifthenelse{\@citet@tag@n > 1}{% multiple tags
  \def\@citet@Tag{Tags}%
 }{%
  \def\@citet@Tag{Tag}%
 }%
 \@citet@Tag~\@citet@tags% return
}
\newcommand{\citet@sep}{, }
\newcommand{\citet@seplast}{ and }
\newcommand{\citet@entry}[2]{\citet@tag{#1}{#2}}
\let\@old@cite\cite
\renewcommand{\cite}[2][]{%
 \def\@citet@detail{\citet@taglist{#1}{\@citet@url}}%
 \ifthenelse{\equal{#2}{sp}}{%
  \def\@citet@url{\citet@url@sp}%
  \def\@citet@bib{\citet@bib@sp}%
 }{\ifthenelse{\equal{#2}{kd}}{%
  \def\@citet@url{\citet@url@kd}%
  \def\@citet@bib{\citet@bib@kd}%
 }{% original cite
  \def\@citet@detail{#1}%
  \def\@citet@bib{#2}%
 }}%
 \ifthenelse{\equal{#1}{}}{%
  \@old@cite{\@citet@bib}%
 }{%
  \@old@cite[\@citet@detail]{\@citet@bib}%
 }%
}
\makeatother

\newcommand{\etale}{{\'etale}}

\newcommand{\DM}{{Deligne--Mumford}}
\newcommand{\Grot}{{Grothendieck}}

\newcommand{\BM}{{Brauer--Manin}}

\newcommand{\TS}{{Tate--Shafarevich}}

\newcommand{\Cart}{{Cartesian}}

\newcommand{\adele}{{ad\`ele}}
\newcommand{\adelic}{{ad\`elic}}

\newcommand{\obs}{\tu{obs}}

\newcommand{\desc}{\tu{desc}}
\newcommand{\conn}{\tu{conn}}
\newcommand{\fin}{\tu{fin}}
\newcommand{\fdesc}{{\fin, \desc}}
\newcommand{\ddesc}{{\desc, \desc}}

\newcommand{\etBr}{{\et, {\Br}}}
\newcommand{\sdesc}{{2\tu{-}\desc}}

\newcommand{\hdesc}[2]{{#1\tu{-}\desc{\ifthenelse{\equal{#2}{}}{}{_#2}}}}

\newcommand{\XA}{X(\bfA_k)}
\newcommand{\YA}{Y(\bfA_k)}
\newcommand{\Xk}{X(k)}

\newcommand{\cXA}{\cX(\cA)}

\newcounter{arrowcounter}

\makeatletter
\zref@newprop{arrownum}{\thearrowcounter}
\zref@addprops{main}{arrownum}

\newcommand{\arrowlabel}[1]{%
    \refstepcounter{arrowcounter}%
      \zref@labelbyprops{#1}{arrownum}%
    \hypertarget{#1}{\tc{\thearrowcounter}}%
          }

\newcommand{\arrowref}[1]{%
            \hyperlink{#1}{\tc{\zref@extract{#1}{arrownum}}}%
              }
\makeatother
\renewcommand{\cXA}{\cX(\bfA_k)}
\newcommand{\finet}{{\fin,\et}}
\newcommand{\dobs}{{\desc,\obs}}
\newcommand{\detBr}{{\desc, \etBr}}

\begin{document}
%title & author
\title[stackcmp]{Comparing local-global obstructions on algebraic stacks}
\author[C. Lv]{Chang Lv}
\address{State Key Laboratory of Cyberspace Security Defense\\
Institute of Information Engineering\\
Chinese Academy of Sciences\\
Beijing 100093, P.R. China}
\email{lvchang@amss.ac.cn}
\author[H. Wu]{Han Wu}
	\address{Hubei Key Laboratory of Applied Mathematics,
	Faculty of Mathematics and Statistics,
	Hubei University,
	No. 368, Friendship Avenue, Wuchang District, Wuhan,
	Hubei 430062, P.R. China}
\email{wuhan90@mail.ustc.edu.cn}
\author[X. Zhang]{Xucheng Zhang}
	\address{School of Mathematical Sciences, Beihang University, Beijing 102206, P.R. China}
\email{xuchengz@buaa.edu.cn}
\date{\today}

\begin{abstract}
Let $k$ be a number field. For quotient stacks $[X/G]$ with $X$ a smooth quasi-projective geometrically
integral $k$-variety and $G$ a linear $k$-group,
 we extend several relations between local-global obstructions previously known for varieties,
  showing that {\etale}-Brauer obstruction is the finest among
  various obstructions (such as iterated descent, finite descent, descent,
  {\BM}).
\end{abstract}
\maketitle

\setcounter{tocdepth}{1}
\tableofcontents

\section{Introduction} \label{intro}
Let $k$ be a number field, and let $X$ be a $k$-variety.
Then one considers the local-global principle for rational points of $X$.
The most developed method is to consider obstructions to it.
More precisely, one constructs various ``nice'' subsets $\XA^{\obs}$ of
 {\adelic} points $\XA$, such that
 $\Xk\subseteq \XA^{\obs}\subseteq \XA$ and hope that
 they capture the local-global principle.

Among others, we are interested in relations between obstructions.
For a smooth, quasi-projective, geometrically integral
 $k$-variety $X$, we have the following relations
 (\cite{harari02groupes,poonen17rational, stoll07finite,
 skorobogatov09descent, demarche09obstruction, poonen10insufficiency,
 hs13homotopy},
 \cite{cdx19comparing, cao20sous}),
\ga{ \label{eq_rel}
\begin{split}
\XA^{\Br} = \XA^{\PGL} =  \XA^\conn = \XA^{\sdesc} = \XA^{h\ZZ} \supseteq \\
   \XA^{h} =  \XA^{\ddesc} =  \XA^{\fdesc} = \XA^{\etBr} = \XA^{\desc}.
\end{split}
}

A natural question is whether these relations hold for algebraic stacks.
The extension of these ideas to algebraic stacks presents both significant
challenges and compelling opportunities. Algebraic stacks provide essential
frameworks for studying moduli problems, quotient constructions, and more
general geometric objects than varieties alone can accommodate. However, this
generalization introduces substantial technical complexities: rational and
{\adelic} points require careful reformulation, classical obstruction theories
must be redeveloped, and new phenomena emerge from stack-theoretic structures.

Recently,  the authors \cite{lh23stackbm} showed that
 $\cXA^{\Br} \subseteq \cXA^{\conn}$ for a
 smooth algebraic $k$-stack
 of finite type that is either {\DM} or
 Zariski-locally the quotient of a smooth geometrically integral
 $k$-variety by a linear $k$-group.
They also showed that $\cXA^\ddesc=\cXA^\desc$ for a quotient stack
 $\cX = [Y/F]$  where  $Y$ is a quasi-projective smooth geometrically integral
 $k$-variety, and $F$ is a finite $k$-group
 \cite{wl24stackdd}.
Nevertheless, a comprehensive comparison of different obstruction theories
on stacks remained an open problem, motivating the present work. The results of this paper establish further relations for a certain class of algebraic stacks.

Another question is which obstructions are the finest (smallest) ones
  on algebraic stacks.
For smooth, quasi-projective geometrically integral
 $k$-variety $X$, \eqref{eq_rel} tells us that $\etBr$ (as well
 as  $h$, $(\ddesc)$, $(\fdesc)$ and $\desc$) is the finest known obstruction.
This paper confirms that $\etBr$ is still the finest known obstruction
 (except $h$ and $h\ZZ$),
 on some class of algebraic
 stacks,
 although we do not know whether $\cXA^{\desc}\subseteq \cXA^{\Br}$ holds
 on these algebraic stacks yet, which is an open question.

% duaring the preparation of this paper, LLLZ showed that ...

This paper provides a systematic development of local-global obstruction
theory for algebraic stacks, demonstrating that under appropriate conditions,
some classical relationships between obstructions extend naturally to the
stack-theoretic setting.
Our main theorem  establishes comprehensive relations:
\thmu{[{Theorem \ref{thm_rel}}]
For $\mathcal{X} = [X/G]$ with $X$ a smooth quasi-projective geometrically
integral $k$-variety and $G$ a linear $k$-group, we have
\[
\begin{tikzcd}[ampersand replacement=\&]
  \cXA^\detBr \arrow[r, equal] \&
  \cXA^\etBr \arrow[r, hookrightarrow] \arrow[d, hookrightarrow] \&
  \cXA^\ddesc \arrow[r, hookrightarrow] \&
  \cXA^\fdesc \arrow[r, equal] \&
  \cXA^\desc \arrow[d, hookrightarrow] \\
  \&
  \cXA^{\Br} \arrow[r, equal] \&
  \cXA^\sdesc \arrow[rr, hookrightarrow] \&
  \&
  \cXA^\conn
\end{tikzcd}
\]
In particular, $\etBr$ is the finest
known obstruction (except  $h$ and $h\ZZ$, whose status is not yet known).
}

Our proofs combine several innovative techniques developed for the
stack-theoretic context. The stable functor framework provides a unifying
language for various obstructions. The stack-theoretic descent formula
 (Proposition \ref{prop_desc}) enables us to construct various composite
  obstructions (such as $\etBr$) as in the case of varieties.
We start by showing the
 descent property of $\cXA^{\Br}$ along an $\SL_n$-torsor
  (Proposition \ref{prop_Br_desc}).
We then show the analogous property of $\cXA^\etBr$ in Proposition \ref{prop_etBr_desc},
 using the existing relations and
 the fact that fundamental groups are preserved under simply
 connected quotients.
Together with functoriality properties of obstructions,
 we obtain the main result (Theorem \ref{thm_rel}).

Our results have immediate implications for arithmetic studies of moduli
 stacks, which naturally appear as quotient stacks $[X/G]$.
Examples include the moduli stack $\mathcal{M}_{g,n}$ of $n$-pointed smooth genus $g$ curves and the moduli stack $\mathrm{Bun}_n^{ss}$ of semistable rank $n$ vector bundles over a curve.
% and $\overline{\mathcal{M}}_g$
% can often be realized as quotient stacks, enabling application of our
% comparison theorems to rational points on moduli spaces.
%
% \item \textbf{Moduli of vector bundles:} For a curve $C$, the moduli stack
% $\mathcal{B}un_{n,d}(C)$ of vector bundles admits quotient descriptions where
% our techniques apply.
%
% \item \textbf{Algebraic cycles and Chow varieties:} Many moduli stacks
% parameterizing algebraic cycles or subschemes can be studied through our
% obstruction framework.
% \end{itemize}
%
In such cases, our results provide tools for determining when the existence of
 local points on a moduli stack implies the existence of global points, with
 potential applications to the arithmetic of moduli problems.

The paper is organized as follows.
Section \ref{ptob} establishes the foundational framework for points and
obstructions on algebraic stacks using stable functors.
Section \ref{quo} discusses when quotients by free actions of algebraic
 groups are schemes, using descent of ample line bundles to show that certain contracted products
 are quasi-projective.
Section \ref{desc}
 develops technical results for quotient stacks, particularly reduction to
 $\SL_n$ quotient stacks and surjectivity properties. Section \ref{cmp} contains our
 main comparison theorems and their proofs.

\section*{Acknowledgments}
The authors would like to thank
 Yang Cao, Ning Guo, Andres Fernandez Herrero, Donghao Li
 and Junchao Shentu for
 helpful discussions.
% and the referees for valuable suggestions.
The work was partially done during the first author's visit to
 the Morningside Center of Mathematics, Chinese
 Academy of Sciences. They thank the Center for its hospitality.

\section{Points and obstructions on algebraic stacks} \label{ptob}
Let us briefly introduce basic notions of rational and {\adelic}
 points, obstructions to local-global principle on algebraic stacks.

We omit set-theoretic issues, as they can be avoided by using
 {\Grot} universes.
Denote by $\Set$ the category of sets, $\Sch$ the category of schemes,
 $\Esp$ the category of algebraic spaces and $\Chp$ the $(2, 1)$-category of
 algebraic stacks.
Here we talk about algebraic spaces and stacks
 in the sense of \cite[025Y,026O]{sp}, which is based on the big fppf
 topology.

For a commutative ring $R$, we also write $R$ for $\Spec R$.
Fix a base scheme $S$ and let $X, T\in \Ob(\Chp_{/S})$ be two
  algebraic stacks over $S$.
The \emph{$T$-points} of $X$ are the isomorphism classes
 of objects in the groupoid $\Hom_{\Chp_{/S}} (T, X)$, denoted by $X(T)$.
Let $k$ be a number field with {\adele} ring $\bfA_k$. The natural inclusion $k\subset \bfA_k$ induces $q: \Spec \bfA_k \ra \Spec k$,
 making $\Spec \bfA_k$ an object of $\Chp_{/k}$.
Let $X\in\Ob(\Chp_{/k})$.
We call $\XA$ (resp. $\Xk$) the \emph{{\adelic} points} (resp.
  \emph{rational points}) of $X$.

\rk{
If $X$ is represented by a scheme, this definition coincides with the classical
 definition.
But in general, when $X$ is an algebraic stack,
 the map  $\Xk\ra \XA$ induced by $q$
 is not necessarily injective.
For example, let $G$ be an affine $k$-group.
Let $*$ be the neutral element in the pointed set
  $BG(\bfA_k) = H_\fppf^1(\bfA_k, G)$.
Then its preimage $\ker(BG(k)\ra BG(\bfA_k))$
 defined by the following {\Cart} diagram
\[
\begin{tikzcd}[ampersand replacement=\&]
  \& \ker(BG(k) \to BG(\mathbf{A}_k))
    \arrow[to=3-1, equal]
    \arrow[r]
    \arrow[d, hook]
  \& \{*\} \arrow[d, hook] \\
  \& BG(k) \arrow[r] \arrow[d, equal]
  \& BG(\mathbf{A}_k) \arrow[d, equal] \\
  \Sha^1(G/k) \arrow[r]
  \& H^1(k, G) \arrow[r]
  \& \check H_\mathrm{fppf}^1(\mathbf{A}_k, G)
\end{tikzcd}
\]
is the {\TS} group $\Sha^1(G/k)$, which is not necessarily trivial.

Nevertheless, by abuse of notation, we also regard $\Xk$ as the image of
 $\Xk\ra\XA$ and write $\Xk\subseteq \XA$.
}

Let us first recall obstructions defined on varieties.
By convention, a \emph{$k$-variety} is a separated $k$-scheme of finite type.
We start from the {\BM} obstruction.
Suppose that  $X$ is a variety over a number field $k$, and
 $\Br X=H_\et^2(X, \bfG_m)$ is the (cohomological) {Brauer--\Grot} group of $X$
 \cite{grothendieck95brauer}.
We have the {\BM} pairing (see, for example, \cite[8.2]{poonen17rational})
\aln{
\XA \tm \Br X & \ra \QQ/\ZZ, \\
((x_v)_v, A) & \mpt \sum_{v\in\Omega_k}\inv_v A(x_v),
}
 where $\bfA_k$ (resp. $\Omega_k$) is the {\adele} (resp. the set of
 all places) of $k$,
  $\inv_v$ is the invariant map at $v$, and
 $A(x_v) = x_v^*A$  is the image of $A$ under the
  pullback functor $x_v^*: \Br X \ra \Br k_v$.
Under this pairing, the {\BM} set (obstruction) is defined to be the subset
 $\XA^{\Br}\subseteq\XA$ orthogonal to the whole $\Br X$, and it plays an important
 role in detecting failures of the local-global principle.

Other obstructions are defined in a similar way.
However, as we mentioned in the introduction, when we want to extend these
 classical obstructions,
 new phenomena emerge from stack-theoretic structures.
To resolve the problem, we first introduce the following notion.
\defi{  [{\cite[2.17]{lv2desc}}] \label{defi_stable}
Let $\sC$ be a $(2, 1)$-category, and $\sD$ be an  ordinary
 category.
Let $F: \sC\ra\sD$ be a functor
 from the underlying ordinary category of $\sC$ to $\sD$.
We say that $F$ is \emph{stable} if
 $F$ is a strict $2$-functor from $\sC$ to $\sD$,
 i.e., for any $1$-morphisms $f$ and $g$ in  $\sC$
  that are $2$-isomorphic, we have $F(f)=F(g)$.
}
\rk{
There is a notion of \emph{stable functor} in the literature
 which means ``having a left adjoint on each slice''.
In this paper, we  keep the
  nomenclature in Definition \ref{defi_stable}.
}

Let $X\in \Chp_{/k}$.
For a stable functor $F: (\Chp_{/k})^\op\ra \Set$ and $\alpha\in F(X)$,
 as in the classical case of Poonen \cite[8.1.1]{poonen17rational},
 one may also define
 the \emph{obstruction given by $\alpha$} to be  the subset
 $\XA^\alpha$ of $\XA$ whose elements are characterized by
\eqn{
\XA^\alpha=\{x\in \XA\mid \alpha(x)\in \im F(q)\},
}
 (which is well-defined by stability of $F$),
 and the  \emph{$F$-set} (or \emph{$F$-obstruction})
 to be the subset $\XA^F$ whose elements are
 characterized by
\eq{ \label{eq_XAF}
\XA^F=\bigcap_{\alpha\in F(X)} \XA^\alpha =
 \{x\in \XA\mid \im F(x)\subseteq \im F(q)\}.
}
Then we have the inclusion
\eqn{
 \Xk\subseteq \XA^F\subseteq  \XA^\alpha\subseteq \XA.
}

\rk{
There is a more functorial way to define $\XA^F$.
Namely, one sees that the map $\Xk\ra \XA$ obviously
 factorizes as
\eq{ \label{eq_obs}
\Xk\ra \Map(F(X), F(k))\tm_{F(q), \Map(F(X), F(\bfA_k)), F}
 \XA\xra{q^*} \XA,
}
 and we define  $\XA^F = \im(q^*)$.
Let $A\in F(X)$. If we replace $F(X)$ in \eqref{eq_obs}
 by $\{A\}$,
 then we define   $\XA^A$ to be the resulting set.
It is easy to verify that these definitions coincide
 with those by \eqref{eq_XAF}.
}

Now let $F = H_\et^n(-, G)$ for some commutative $k$-group $G$.
Then $F: (\Chp_{/k})^\op \ra \Set$
 is stable (cf. \cite[2.29]{lv2desc})
 and the
 resulting $F$-set  $\XA^{H_\et^n(-, G)}$ coincides with the classical definition
 if  $X$ is a $k$-variety.
In particular, $F = \Br = H_\et^2(-, \bfG_m)$ gives
 the \emph{{\BM} obstruction} of $\XA^{\Br}$.

Next let $G$ be a $k$-group and  $T\in\Ob(\Chp_{/k})$, write $G_T=G\tm_k T$.
Then one considers fppf
 $G_T$-torsors
 over $T$
 (see, for example, \cite[Sec. 2.3]{lh23stackbm}).
Let $\Tors(X, G)$ be the groupoid of  fppf $G$-torsors over $X$.
There is an isomorphism of pointed sets (cf.
 \cite[III. 3.6.5 (5), IV. 3.4.2 (i)]{giraud71cohnonab})
\eqn{
  \check H_\fppf^1(X, G)\xra{\sim} \Tors(X, G)/\cong.
}
It turns out that
 $F=\check H_\fppf^1(-, G): (\Chp_{/k})^\op\ra \Set$ is also a stable functor
 (cf. \cite[3.14]{lv2desc})
 and the
 resulting $F$-set  $\XA^{\check H_\fppf^1(-, G)}$ coincides with the classical
 definition
 if
 $X$ is a $k$-variety. See \cite[3.17 (i)]{lv2desc}.
Define  the \emph{descent obstruction} of $X$ to be
\eqn{
\XA^\desc=\bigcap_{\text{linear $k$-group $G$}} \XA^{\check H_\fppf^1(-, G)}.
}
Similarly, define
\eqn{
\XA^\conn=\bigcap_{\text{connected linear $k$-group $G$}}
 \XA^{\check H_\fppf^1(-, G)},
}
 and the \emph{second descent obstruction}
   \cite[4.2]{lv2desc}
\eqn{
\XA^{\sdesc}=
 \bigcap_{\text{commutative linear $k$-group $G$}}
 \XA^{H_\et^2(-, G)}.
}

Let $f: Y\xra{G} X$ be a $G$-torsor over $X$
  (which is also in $\Chp_{/k}$ by \cite[Lem. 2.10]{lh23stackbm}),
  and denote also by $f$
  the element of
  $\check H_\fppf^1(X, G)$ corresponding to the class of
  $Y\xra{G} X$.
\prop{ \label{prop_desc}
We have
\eqn{
\XA^f=\bigcup_{\s\in H^1(k, G)}
 f^\s(Y^\s(\bfA_k)),
}
 where $f^\s: Y^\s\xra{G^\s} X$  is the twist of
  $f$ $($see \cite[3.30]{lv2desc}$)$.
}
\pf{
For $X$ being a  $k$-variety,  this is well-known.
For algebraic stacks, this follows from general descent by torsors
 \cite[3.20]{lv2desc}.
}
From this we may define  the \emph{{\etale} Brauer  obstruction} of $X$, to be
\eqn{
\XA^\etBr=\bigcap_{f:\ Y\xra{F} X \text{ torsor under
 finite $k$-group $F$}} \bigcup_{\s\in H^1(k, F)}
 f^\s(Y^\s(\bfA_k)^{\Br}),
}
 the \emph{finite  descent  obstruction}
\eqn{
\XA^\fdesc=\bigcap_{f:\ Y\xra{F} X \text{ torsor under
 finite $k$-group $F$}} \bigcup_{\s\in H^1(k, F)}
 f^\s(Y^\s(\bfA_k)^{\desc}),
}
 and the \emph{iterated descent  obstruction} of $X$, to be
\eqn{
\XA^\ddesc=\bigcap_{f:\ Y\xra{G} X \text{ torsor under
 linear $k$-group $G$}} \bigcup_{\s\in H^1(k, G)}
 f^\s(Y^\s(\bfA_k)^{\desc}).
}

All of these definitions of obstructions coincide with  the
 classical ones in the case that $X$ is a $k$-variety.

\defi{ \label{defi_obs_functorial}
Suppose that attached to each $X\in \Chp_{/k}$,
 we have a subset $\XA^\obs\subseteq \XA$.
We say that $\obs$ is \emph{functorial} on $\Chp_{/k}$
 if the assignment $X\mpt \XA^\obs$ gives a functor
\eqn{
-(\bfA_k)^\obs: \Chp_{/k}\ra \Set.
}}

\rk{ \label{rk_functorial}
In particular, if
 $\obs: (\Chp_{/k})^\op\ra \Set$
 is a stable functor, and for $X\in \Chp_{/k}$, let
 $\XA^\obs$ be the corresponding
 obstruction set,
then  $\obs$ is functorial on  $\Chp_{/k}$,
}

\lemm{ \label{lemm_func}
Let $\obs$ be a functorial obstruction  on $\Chp_{/k}$.
Define
\eqn{
\XA^\dobs=\bigcap_{f:\ Y \xra{G} X \mathrm{\ torsor\ under\ linear\ } k\text{-}\mathrm{group\ } G} \bigcup_{\s\in H^1(k, G)}
 f^\s(Y^\s(\bfA_k)^{\obs}).
}
Then we have
\enmt{[\upshape (i)]
\item \label{it_func} $(\dobs)$ is also a functorial obstruction
 on $\Chp_{/k}'$, and
\item \label{it_cap} $\XA^\dobs\subseteq \XA^\desc\cap \XA^\obs$.
}}
\pf{
There is a more general statement and proof in \cite[Thm. 5.10]{lv2desc}.
For the reader's convenience, we give a specific proof here.

For \eqref{it_func},
 let $f: X'\ra X$ be a  $1$-morphism in $\Chp_{/k}$,
 and we need to show that $f(X'(\bfA_k)^\dobs)\subseteq \XA^\dobs$.
Let $G$ be a linear $k$-group and $\s\in H^1(k, G)$.
Let $g: Y\xra{G} X$ be a $G$-torsor, and $g^\s: Y^\s\xra{G^\s} X$ be
 its twist by $\s$.
Let $Y'=Y\tm_X X'$.
Then we have the pullback square
\[
\begin{tikzcd}[ampersand replacement=\&]
  Y'^\s \arrow[r, "g'^\s" above] \arrow[d, "f'^\s" left]
  \& X' \arrow[d, "f" right] \\
  Y^\s \arrow[r, "g^\s" above]
  \& X
\end{tikzcd}
\]
 with $g'^\s: Y'^\s\ra X'$ the twist of the $G$-torsor $g':Y'\ra X$ by $\s$.
Since $\obs$ is functorial on $\Chp_{/k}$, we have
 $f'(Y'^\s(\bfA_k)^\obs)\subseteq Y^\s(\bfA_k)^\obs$.
Let $G$ range over all linear $k$-groups and $g$ over all $G$-torsors on $X$, we obtain that
\aln{
f(X'(\bfA_k)^\dobs) &\subseteq f\left( \bigcap_{G,\ g:\ Y\xra{G} X}
 \bigcup_{\s\in H^1(k, G)} g'^\s(Y'^\s(\bfA_k)^\obs) \right) \\
 &\subseteq \bigcap_{G,\ g:\ Y\xra{G} X}
  \bigcup_{\s\in H^1(k, G)} fg'^\s(Y'^\s(\bfA_k)^\obs)  \\
 &= \bigcap_{G,\ g:\ Y\xra{G} X}
  \bigcup_{\s\in H^1(k, G)} g^\s f'^\s(Y'^\s(\bfA_k)^\obs)  \\
 &\subseteq \bigcap_{G,\ g:\ Y\xra{G} X}
  \bigcup_{\s\in H^1(k, G)} g^\s(Y^\s(\bfA_k)^\obs)  \\
 &= X(\bfA_k)^\dobs,
}
 which is desired.

For \eqref{it_cap}, the inclusion $\XA^\dobs\subseteq \XA^\desc$ is due to
 Proposition \ref{prop_desc} and
 $\XA^\dobs\subseteq \XA^\obs$ is due to the functoriality of $\obs$.

The proof is complete.
}

\section{Schematic quotients}  \label{quo}
In this paper, we need to investigate some objects that are
 quotients by algebraic groups, and these are not necessarily schemes,
 even when a reductive group acts freely on a smooth quasi-projective variety
 (see \cite[Cor. 6]{Kollar-non-scheme}).
This section is devoted to situations where quotients under free actions, a
 priori merely algebraic spaces
 (see \cite[071S]{sp}), are schemes. Typically, there are two ways to prove that an algebraic space is a scheme:
\begin{itemize}
\item
either equip it with an ample line bundle (any quasi-projective algebraic space is a scheme, see \cite[Prop. 5.5.27]{alper_stacks_moduli}), cf. Theorem \ref{thm_normal-qs_quotient}, or
\item
show it admits a Zariski open covering by schemes (by definition), cf. Lemma \ref{lem_split_scheme}.
\end{itemize}

For the first approach, the key technical result is:
\prop{[\cite{mumford-git}, Chap. 7, Prop. 7.1]\label{thm_quotient_scheme}
Let $S$ be a locally noetherian scheme. Let $G$ be a group scheme, flat and of finite type over $S$. Let $X$ and $Y$ be $G$-schemes of finite type over $S$, and let $f: X \to Y$ be a $G$-equivariant $S$-morphism. Assume there exist
\begin{itemize}
\item
a $G$-torsor $\varphi_Y: Y \to Q$ over a quasi-projective $S$-scheme $Q$.
\item
a $G$-equivariant line bundle $L$ on $X$ which is relatively ample for $f$.
\end{itemize}
Then there exists a $G$-torsor $\varphi_X: X \to P$ over a quasi-projective $S$-scheme $P$. More precisely
\begin{itemize}
\item
the $G$-equivariant $S$-morphism $f: X \to Y$ descends to an $S$-morphism $\bar{f}: P \to Q$ and the resulting commutative diagram
\[
\begin{tikzcd}[ampersand replacement=\&]
X \ar[r,"f"] \ar[d,dashed,"\varphi_X"'] \& Y \ar[d,"\varphi_Y"] \\
P \ar[r,dashed,"\bar{f}"'] \& Q \ar[ul,phantom,"\lrcorner"]
\end{tikzcd}
\]
is Cartesian.
\item
the $G$-equivariant line bundle $L$ on $X$ descends to a line bundle $\bar{L}$ on $P$ which is relatively ample for $\bar{f}$.
\end{itemize}
}

As a corollary, we obtain
\cor{\label{cor_torsor_quot}
Let $S$ be a locally noetherian scheme. Let $G$ be a group scheme, flat and of finite type over $S$. Let $X$ and $Y$ be $G$-schemes of finite type over $S$. Let $G$ act diagonally on $Y \times_S X$. Assume there exist
\begin{itemize}
\item
a $G$-torsor $\varphi_Y: Y \to Q$ over a quasi-projective $S$-scheme $Q$.
\item
a $G$-equivariant $S$-ample line bundle $L$ on $X$.
\end{itemize}
Then there exists a $G$-torsor $Y \times_S X \to (Y \times_S X)/G$ over a quasi-projective $S$-scheme $(Y \times_S X)/G$.
}
\rk{
The fppf quotient $(Y \times_S X)/G$ of $Y \times_S X$ by the diagonal $G$-action is called the \emph{contracted product} in this paper, denoted by $Y \times_S^G X$. This result generalizes \cite[Lem. 2.2.3]{torsor} or \cite[6.5.6.3]{poonen17rational}, where $X$ is assumed to be
 \emph{affine} over $S$.
Indeed, for any $G$-scheme $X$, the structure sheaf $\mathcal{O}_X$ is always
 $G$-equivariant and it is $S$-ample if and only if $X$ is a quasi-affine
 $S$-scheme (\cite[0891]{sp}).
 % namely, we use GIT instead of {\Grot} fppf descent to remove the affine
 % condition.
}
\begin{proof}
Let $p_2: Y \times_S X \to X$ be the second projection; it is $G$-equivariant. For the $G$-equivariant $S$-ample line bundle $L$ on $X$, the pullback $p_2^*L$ is a $G$-equivariant line bundle on $Y \times_S X$ which is relatively ample for the first projection $p_1: Y \times_S X \to Y$.
Then Proposition \ref{thm_quotient_scheme} completes the diagram to a Cartesian one
\[
\begin{tikzcd}
Y \times_S X \ar[r,"p_1"] \ar[d,dashed] & Y \ar[d,"\varphi_Y"] \\
(Y \times_S X)/G \ar[r,dashed] & Q \ar[ul,phantom,"\lrcorner"]
\end{tikzcd}
\]
and $Y \times_S X \to (Y \times_S X)/G$ is a $G$-torsor over a quasi-projective $S$-scheme $(Y \times_S X)/G$.
\end{proof}
In the application of Corollary \ref{cor_torsor_quot}, we take the torsor $Y \to Q$ to be the fppf quotient $G \to G/H$. Therefore, the fppf quotient $G/H$ should be $S$-quasi-projective. However, this fails even when $S$ is a 2-dimensional affine space, as shown by the counterexample in \cite[X 13]{Raynaud-quasi-projective}. For bases of dimension at most one, the representability of $G/H$ is established by S. Anantharaman, and under additional assumptions it is known to be quasi-projective. %this is precisely why we restrict to bases of dimension at most one.
% this is still an open problem over an arbitrary base scheme $S$ (raised by M. Raynaud, see \cite[XV, 2.ii) and 6]{Raynaud-quasi-projective}). Nevertheless, it is known to hold in several cases. \cite[Chapitre IV, 4.C. Th'{e}or`{e}m]{Anantharaman-dim-1}Nevertheless, we know it holds in the following case.
\lemm{[Quasi-projectivity of $G/H$]\label{lem_quasi-projective-G/H}
Let $S$ be a locally noetherian semilocal Pr\"{u}fer scheme (i.e. $S=\Spec(R)$ for some semilocal Dedekind domain $R$; in particular $\dim S \leq 1$). Let $G$ be a smooth group scheme of finite type over $S$. Let $H \subseteq G$ be a closed subgroup scheme, flat over $S$. Then the fppf quotient $G/H$ is a quasi-projective $S$-scheme.
}
\pf{
By \cite[Ch. IV, 4.C. Th]{Anantharaman-dim-1}, the fppf quotient $G/H$ is an $S$-scheme; it is also $S$-separated because $H \subseteq G$ is closed. Indeed, after the fppf base change $G \times_S G \to (G/H) \times_S (G/H)$, the diagonal of $G/H$ becomes the closed immersion $G \times_S H \hookrightarrow G \times_S G$ given by $(g,h) \mapsto (g,gh)$. Under our assumptions on the base $S$, we know the $S$-separated $G$-homogeneous space $G/H$ is quasi-projective by \cite[Cor. 6.12]{Guo-quasi-proj}. %\cite[Ch. VI, Cor. 2.5]{Raynaud-quasi-projective}
}
%\cite[Corollary 11.5]{Pappas-Zhu-2-dim} $S$ is an affine excellent regular scheme of $\dim S=2$, and $G$ and $H$ are smooth affine group schemes over $S$ with connected fibers.
In the application of Corollary \ref{cor_torsor_quot}, we also need a $G$-equivariant ample line bundle; its existence is guaranteed essentially by Sumihiro's theorem.
%\lemm{[\cite{Sumihiro-II}, Thm. 1.6]\label{lem_sumihiro-II}
%Let $S$ be a scheme. Let $G$ be a smooth affine group scheme over $S$ with connected fibers (e.g. $G=\mathrm{GL}_{n,S}$), and $X$ be a normal noetherian $G$-scheme over $S$. Then for any line bundle on $X$, some positive power is $G$-equivariant. In particular, any noetherian normal $G$-scheme quasi-projective over $S$ admits a $G$-equivariant $S$-ample line bundle.
%}
\lemm{\label{lem_sumihiro-II}
Let $S$ be a locally noetherian scheme. Let $G$ be a smooth affine group scheme of finite type over $S$, and let $G^\circ \subseteq G$ be the relative identity component. Let $X$ be a noetherian $G$-scheme over $S$. Suppose
\begin{itemize}
\item
the fppf quotient $G/G^\circ$ is a finite locally free $S$-scheme (e.g. if $S=\mathrm{Spec}(k)$, in which case the fppf quotient $G/G^\circ$ is finite \'{e}tale over $k$).
\item
the scheme $X$ admits a $G^\circ$-equivariant (resp. $S$-ample) line bundle (e.g. if $X$ is normal, then every line bundle on $X$ has a positive power that is $G^\circ$-equivariant, see \cite[Thm. 1.6]{Sumihiro-II}).
\end{itemize}
Then $X$ admits a $G$-equivariant (resp. $S$-ample) line bundle. %In particular, any noetherian normal $G$-scheme quasi-projective over $S$ admits a $G$-equivariant $S$-ample line bundle.
}
\pf{
Let $M$ be a $G^\circ$-equivariant line bundle on $X$. Let $p_2: G \times_S X \to X$ be the second projection. The line bundle $p_2^*M$ on $G \times_S X$ descends along the $G^\circ$-torsor $G \times_S X \to G \times_S^{G^\circ} X$, say descends to a line bundle $M^\circ$ on $G \times_S^{G^\circ} X$. Note that there is a canonical isomorphism
\[
G \times_S^{G^\circ} X \xrightarrow{\sim} (G/G^\circ) \times_S X \text{ mapping } [g,x] \mapsto (gG^\circ,g.x).
\]
Thus we can view $M^\circ$ as a line bundle on $(G/G^\circ) \times_S X$. Let $q: (G/G^\circ) \times_S X \to X$ be the second projection; it is finite locally free since $G/G^\circ \to S$ is. Define the norm line bundle on $X$
\[
N:=\mathrm{Nm}_{q}(M^\circ):=\det\left(q_*M^\circ\right) \otimes \det\left(q_*\mathcal{O}_{(G/G^\circ) \times_S X}\right)^{-1}.
\]
Because all operations involved are functorial for equivariant finite locally free morphisms, the line bundle $N$ is $G$-equivariant.

Suppose $M$ is $S$-ample; we show $N$ is also $S$-ample, which can be checked after an fppf base change on $S$. After such a base change, we may assume that the fppf quotient $G/G^\circ$ is the disjoint union of $d$ copies of $S$, and these components have lifts $g_1,\dots,g_d$ in $G$. Then
\[
G \times_S^{G^\circ} X \cong (G/G^\circ) \times_S X \cong \coprod_{i=1}^d X
\]
and $M^\circ|_{X_i} \cong g_i^*M$, where $g_i: X \to X$ is an $S$-automorphism. Consequently,
\[
N \cong \bigotimes_{i=1}^d g_i^*M.
\]
Each $g_i^*M$ is $S$-ample, and since a finite tensor product of relatively ample line bundles is relatively ample, we conclude that $N$ is $S$-ample.
}
\rk{
Even when $G/G^\circ$ is finite \'{e}tale, one cannot claim that every $G^\circ$-equivariant line bundle $M$ has a positive power that is $G$-equivariant. For example, over a field, take
\[
G=\mathbf{Z}/2 \ \curvearrowright \ \mathbf{P}^1 \times \mathbf{P}^1=X,
\]
with $G$ exchanging the factors, and
\[
M=\mathcal{O}_X(a,b) \text{ with } a \neq b.
\]
Then no positive power $M^{\otimes n}$ is $G$-invariant. The norm instead produces
\[
\mathrm{Nm}(M)=M \otimes \sigma^*M \cong \mathcal{O}_X(a+b,a+b),
\]
which is $G$-equivariant.
}
Finally, we list some properties that descend along $G$-torsors (or more generally, along fppf coverings).
\lemm{\label{lem_desend_prop}
Let $S$ be a scheme. Let $X \to Y$ be an fppf covering of $S$-schemes (e.g. any $G$-torsor for a group scheme $G$ flat and locally of finite presentation over $S$). If $X$ is normal, then $Y$ is normal. Moreover, if $X \to S$ is smooth (resp. locally of finite type/presentation, has geometrically integral fibers), so is $Y \to S$.
}
\pf{
Normality: For any $y \in Y$, choose $x \in X$ over $y$. The induced morphism $\cO_{Y,y} \to \cO_{X,x}$ is flat and local, hence faithfully flat. Since $X$ is normal, the local ring $\cO_{X,x}$ is normal and so is $\cO_{Y,y}$. This shows that $Y$ is normal. Smoothness: \cite[05B5]{sp}. Locally of finite type/presentation: \cite[02KL]{sp}.

Having geometrically integral fibers: note that a morphism has geometrically integral fibers if and only if its geometric fibers are integral; we directly reduce to the case that $S=\mathrm{Spec}(k)$ for some algebraically closed field $k$. Then the integrality of $X$ implies that of $Y$ (see \cite[06QM]{sp}).
}
Now we are in a position to state the main result in this section, which may be viewed as a relative version of \cite[Lem. 2.11]{Brion-action}.
\thm{\label{thm_normal-qs_quotient}
Let $S$ be a locally noetherian semilocal Pr\"{u}fer scheme. Let $G$ be a smooth group scheme of finite type over $S$. Let $H \subseteq G$ be a closed subgroup scheme, smooth and affine over $S$ such that the fppf quotient $H/H^\circ$ is a finite locally free $S$-scheme. Let $X$ be a normal $H$-scheme quasi-projective over $S$. Let $H$ act diagonally on $G \times_S X$. Then the quotient
\[
G \times_S^H X:=(G \times_S X)/H
\]
is a normal quasi-projective $S$-scheme. Moreover, if $S=\mathrm{Spec}(k)$, $G$ is connected and $X$ is a smooth quasi-projective geometrically integral $k$-variety, so is $G \tm_k^H X$.
}
\begin{proof}
We verify the hypotheses of Corollary \ref{cor_torsor_quot} for the diagram
\[
\begin{tikzcd}[ampersand replacement=\&]
G \times_S X \ar[r,"p_1"] \& G \ar[d,"H\text{-tor}"] \\
 \& G/H
\end{tikzcd}
\]
The fppf quotient $G \to G/H$ is an $H$-torsor and the fppf quotient $G/H$ is a quasi-projective $S$-scheme by Lemma \ref{lem_quasi-projective-G/H}. By Lemma \ref{lem_sumihiro-II} there exists an $H$-equivariant $S$-ample line bundle on $X$. By Corollary \ref{cor_torsor_quot} there exists an $H$-torsor $G \times_S X \to G \times_S^H X$ over a quasi-projective $S$-scheme $G \times_S^H X$. Since $G \to S$ is smooth, the second projection $G \times_S X \to X$ is smooth with normal target $X$, so $G \times_S X$ is normal and we can apply Lemma \ref{lem_desend_prop} to conclude the normality of $G \times_S^H X$.

Suppose $S=\mathrm{Spec}(k)$ and $G$ is connected. If $X$ is a smooth geometrically integral $k$-variety, so is $G \times_k X$. By Lemma \ref{lem_desend_prop}, the fppf quotient $G \tm_k^H X$ is smooth and geometrically integral.
\end{proof}
%\rk{
%It is clear from the proof that the condition ``$X$ is a normal $H$-scheme quasi-projective over $S$'' in Theorem \ref{thm_normal-qs_quotient} can be weakened to ``$X$ is an $H$-scheme of finite type over $S$ with an $H$-equivariant $S$-ample line bundle''.
%}

\cor{\label{cor_contracted_rep}
Let $S$ be a locally noetherian scheme. Let $H$ be a smooth affine group scheme over $S$. Let $X$ be a normal $H$-scheme quasi-projective over $S$. Let $f: X \to [X/H]$ be the $H$-torsor.
\enmt{[\upshape (i)]
\item
\label{it_twist_rep} For any $\s \in H^1(S,H)$, the twist $X^\s$ is a normal
 quasi-projective $S$-scheme, and the projection $f^\s: X^\s \ra [X/H]$ is an
 $H^\s$-torsor.

 Moreover, if $S=\mathrm{Spec}(k)$ and $X$ is a smooth quasi-projective geometrically integral $k$-variety, so is $X^\s$. % FIXME: geo int condition ???  SOLVED
\item
\label{it_H_rep} Suppose $S$ is a semilocal Pr\"{u}fer scheme and the fppf quotient $H/H^\circ$ is a finite locally free $S$-scheme. Let $G$ be a smooth group scheme of finite type over $S$, and suppose that it contains $H$ as a closed subgroup scheme. The contracted product $G \tm_S^H X$ is a normal
 quasi-projective $S$-scheme, and the projection $G \tm_S^H X \to [X/H]$
 is a $G$-torsor.
 % FIXME: geo int condition ???  SOLVED $G$ is connected
}}
\pf{
For \eqref{it_twist_rep}, let $Y \to S$ be an $H$-torsor representing the class $\s \in H^1(S,H)$. Consider the following Cartesian diagram
\[
\begin{tikzcd}[ampersand replacement=\&]
Y \times_S X \ar[r] \ar[d] \& X \ar[d,"f"] \\
Y \times_S {[}X/H{]} \ar[r] \ar[d] \& {[}X/H{]} \ar[d] \ar[ul,phantom,"\lrcorner"] \\
Y \ar[r,"\s"'] \& S \ar[ul,phantom,"\lrcorner"]
\end{tikzcd}
\]
By definition, $X^\s:=Y \times_S^H X$ and $H^\s:= Y \times_S^H H$; the latter is the inner twist of $H$ by $\s$, with $H$ acting on itself by conjugation. Applying Corollary \ref{cor_torsor_quot} to the diagram
\[
\begin{tikzcd}[ampersand replacement=\&]
Y \times_S X \ar[r,"p_1"] \& Y \ar[d,"H\text{-tor}"] \\
 \& S
\end{tikzcd}
\]
gives an $H$-torsor $Y \times_S X \to Y \times_S^H X=X^\s$ over a quasi-projective $S$-scheme $X^\s$. The second projection $Y \times_S X \to X$ is smooth (as the base change of $Y \to S$) with normal base $X$, so $Y \times_S X$ is normal and hence $X^\s$ is normal by Lemma \ref{lem_desend_prop}. Moreover, the second projection $Y \times_S X \to X$ is $H$-equivariant so it descends to the fppf quotients $f^\s: X^\s \to [X/H]$. The twisted group scheme $H^\s=Y \times_S^H H \cong \mathrm{Aut}_H(Y)$ acts on $X^\s= Y \times^H_S X$ via its natural action on the first factor. This action makes $f^\s : X^\s \to [X/H]$ an $H^\s$-torsor. Indeed, the $H$-torsor $f: X \to [X/H]$ is the pullback of the universal $H$-torsor over $BH$ along $[X/H] \to BH$, and twisting by $\s$ replaces the structure group $H$ by its inner form $H^\s$.

Since $Y$ is smooth over $S$, if $X$ is smooth over $S$, so is $Y \times_S X$. By Lemma \ref{lem_desend_prop}, the twist $X^\s$ is smooth over $S$. If $S=\Spec(k)$, after base change to $\bar{k}$, the $H$-torsor $Y$ becomes trivial, and hence
\[
(X^\s)_{\bar k}=Y_{\bar k} \times^{H_{\bar k}}_{\bar k} X_{\bar k}=X_{\bar k}.
\]
Therefore $X^\sigma$ is geometrically integral whenever $X$ is geometrically integral.

For \eqref{it_H_rep}, the contracted product $G \tm_S^H X$ is a normal quasi-projective $S$-scheme by Theorem \ref{thm_normal-qs_quotient}. The $H$-equivariant second projection $G \tm_S X \to X$ descends to quotients $G \tm_S^H X \to [X/H]$, which is a $G$-torsor as it is the extension of the structure group of the $H$-torsor $f: X \to [X/H]$ along the subgroup $H \subseteq G$.
}
\lemm{\label{lemm_quo_SLn}
Let $S$ be a locally noetherian semilocal Pr\"{u}fer scheme. Let $\cX=[X/G]$ be a quotient stack over $S$, where $X$ is an algebraic space over $S$, and $G$ is a smooth affine group scheme over $S$ such that the fppf quotient $G/G^\circ$ is a finite locally free $S$-scheme. Then there exist an integer $n>0$ and an algebraic space $Y$ over $S$ such that $\cX \cong [Y/\SL_{n,S}]$. If $X$ is a normal quasi-projective $S$-scheme, then $Y$ can be arranged to be so; the same conclusion holds if $S=\mathrm{Spec}(k)$ and $X$ is a smooth quasi-projective geometrically integral $k$-variety.
}
\pf{
The assumptions on $S$ and $G$ guarantee that there is a closed immersion $G \hookrightarrow \GL_{n-1,S}$ for some integer $n>1$ (see \cite[Exp. $\mathrm{VI_B}$, 11.11.1]{SGA3} or \cite[Cor. 3.2]{Thomason-embedding}), and composing with $\GL_{n-1,S} \hookrightarrow \SL_{n,S}$ yields an embedding into $\SL_{n,S}$. Set
\[
Y:=\SL_{n,S} \tm_S^G X :=(\SL_{n,S} \tm_S X)/G.
\]
The $G$-action on $\SL_{n,S} \tm_S X$ is free, so $Y$ is an algebraic space over $S$, and
\[
[X/G] \cong \left[\left((\SL_{n,S} \tm_S X)/G\right)/\SL_{n,S}\right] = [Y/\SL_{n,S}].
\]
The remaining assertion follows from Theorem \ref{thm_normal-qs_quotient}.
}
To close this section, we record another situation where the fppf quotient $G \times_S^H X$ is a scheme, although it will not be used in this paper.
\lemm{\label{lem_split_scheme}
Let $S$ be a locally noetherian semilocal Pr\"{u}fer scheme. Let $G$ be a smooth group scheme of finite type over $S$. Let $H \subseteq G$ be a closed subgroup scheme, flat over $S$. Let $X$ be an $H$-scheme over $S$. If the $H$-torsor $\pi: G \to G/H$ is Zariski locally trivial (e.g. if $H$ is special), then $G \times_S^H X$ is an $S$-scheme.
}
\pf{
Note that we have a commutative diagram
\[
\begin{tikzcd}[ampersand replacement=\&]
G \times_S X \ar[r,"p_1"] \ar[d,"\varphi"'] \& G \ar[d,"\pi"] \\
(G \times_S X)/H \ar[r] \& G/H
\end{tikzcd}
\]
and $\varphi: G \times_S X \to (G \times_S X)/H$ is an $H$-torsor over an algebraic space $(G \times_S X)/H$. Once again, the fppf quotient $G/H$ is a quasi-projective $S$-scheme by Lemma \ref{lem_quasi-projective-G/H}. Choose a Zariski open covering $\{U_i \to G/H\}$ that trivializes the $H$-torsor $\pi: G \to G/H$, i.e., each $U_i \hookrightarrow G/H$ is an open immersion of $S$-schemes and
\[
\pi|_{U_i}: \pi^{-1}(U_i) \cong U_i \times_S H \to U_i
\]
is the trivial $H$-torsor over $U_i$. Since $\varphi$ is flat, locally of finite presentation (hence open) and surjective, the images
\[
\varphi(p_1^{-1}(\pi^{-1}(U_i))) \hookrightarrow (G \times_S X)/H
\]
form a Zariski open covering of $(G \times_S X)/H$. To conclude, we show each $\varphi(p_1^{-1}(\pi^{-1}(U_i)))$ is an $S$-scheme. Indeed, we have
\[
\varphi(p_1^{-1}(\pi^{-1}(U_i))) = p_1^{-1}(\pi^{-1}(U_i))/H \cong (U_i \times_S H \times_S X)/H \xrightarrow{\sim} U_i \times_S X,
\]
where the last morphism is $[(u,h,x)] \mapsto (u,h.x)$, with the inverse $(u,x) \mapsto [(u,e_H,x)]$. The fiber product $U_i \times_S X$ is an $S$-scheme as both $U_i$ and $X$ are.
}

\section{Descent for \texop{$\Br$}{Br} and
 \texop{$\etBr$}{\et, Br} using \texop{$\SL_n$}{SLn}} \label{desc}

\lemm{\label{lemm_H1_SLn_field}
We have $H^1(L, \SL_n) = \{*\}$ for any field $L$.
}
\pf{
Consider the exact sequence of $L$-groups
\[
1 \to \SL_n \to \GL_n \xra{\det} \bfG_m \to 1.
\]
Taking Galois cohomology yields
\[
\GL_n(L) \xra{\det} L^\times \to H^1(L, \SL_n) \to H^1(L, \GL_n).
\]
By Hilbert's Theorem 90 \cite[Lem.~1, p.~129]{serre94cohomologie},
 we have $H^1(L, \GL_n) = \{*\}$.
The map $\GL_n(L) \to L^\times$ is surjective because any $a \in L^\times$ is the
determinant of $\diag(a,1,\dots,1) \in \GL_n(L)$. Hence $H^1(L, \SL_n)$ is trivial.
}

\lemm{ \label{lemm_H1_SLn}
We have
 $\check H_\fppf^1(\bfA_k, \SL_n) = \{*\}$.
}
\pf{
Since the natural  map
\eqn{
\check H_\fppf^1(\bfA_k, \SL_n) \ra \bigoplus_v H^1(k_v, \SL_n)
}
 is injective (a consequence of
 \cite[Thm. 2.18]{cesnavicius15poitou}), by Lemma \ref{lemm_H1_SLn_field},
the Galois cohomology $H^1(k_v, \SL_n)$ is trivial, and this lemma follows.
}

\prop{ \label{prop_Ak_desc}
Let $\cX=[Y/\SL_n]$ where $Y$ is a $k$-scheme.
Then the natural map  $Y\ra \cX$ induces a surjective map
 $\YA\ra \cXA$.
}
\pf{We take an {\adelic} point $x\in \cXA$. Its fiber under the morphism
 $Y\ra \cX$ is an $\SL_n$-torsor over $\bfA_k$. By Lemma \ref{lemm_H1_SLn}, it is a
 trivial $\SL_n$-torsor. So there exists an {\adelic} point of $\YA$ lying over $x$,
 which implies our proposition.
}

In the following lemma, we  generalize  \cite[Lem. 2.1]{cdx19comparing} to the quotient stacks.
\lemm{ \label{lemm_Br_iso}
	Let $H$ be a semi-simple, simply connected group over $k$. Let $X$ be a smooth geometrically integral $k$-variety with an $H$-action $\rho\colon H\times X\to X$. Let $[X/H]$ be the quotient stack.
Then the natural projection $f\colon X\to [X/H]$
% be $p_2: \SL_n\tm_k Y\ra Y$
 induces an isomorphism $f^*\colon \Br[X/H]\to\Br X.$
% $p_2^*:\Br Y\xra{\sim} \Br(\SL_n\tm_k Y)$.
}

\pf{	  We use Sansuc's exact sequence for algebraic stacks
	(\cite[Thm. 3.1]{lh23stackbm})
	for the torsor $f: X\xra{H} [X/H]$.
	Then we obtain an exact sequence
	\eqn{
		0 =  \Pic H\ra \Br[X/H]\xra{f^*} \Br X\xra{\rho^*-p_2^*}
		\Br(H\tm X).
	}
	Consider morphisms
	\[
	\begin{tikzcd}[ampersand replacement=\&]
  X \arrow[r, hookrightarrow, "i" above] \&
  H\tm_k X \arrow[r, "p_2" above, "\rho" below] \&
  X
\end{tikzcd}
	\]where the first injection is induced by the unit of $H$.
	Then we have a commutative diagram
	\[
	\begin{tikzcd}[ampersand replacement=\&]
  \& X \arrow[dr, equal] \arrow[dl, hookrightarrow, "i" above] \& \\
  H\tm_k X \arrow[rr, "p_2" above, "\rho" below] \& \& X
\end{tikzcd}
	\] which induces the following diagram
	\[
	\begin{tikzcd}[ampersand replacement=\&]
  \& \Br X \arrow[dr, equal] \& \\
  \Br (H\tm_k X) \arrow[ur, "i^*"] \& \& \Br X \arrow[ll, "p_2^*" above, "\rho^*" below]
\end{tikzcd}
	\]
	By \cite[Lem. 2.1]{cdx19comparing}, $p_2^*$ is an isomorphism,
  hence $i^*$ also an isomorphism,
	and $p_2^*=\rho^*$. Hence, $f^*:\Br[X/H]\to \Br X$ is an isomorphism.
}

\prop{ \label{prop_Br_desc}
Let $\cX=[Y/\SL_n]$ where $Y$ is  a smooth geometrically integral
 $k$-variety.
Let  $f$ be the natural map  $Y\ra \cX$.
Then for any $x\in \cXA^{\Br}$, any lift $y\in \YA$ along $f$
 must lie in $\YA^{\Br}$.
In particular, $f$ induces a surjective map
 $\YA^{\Br}\ra \cXA^{\Br}$.
}
\pf{
By the functoriality of $\Br$ (see Remark \ref{rk_functorial}),
 there is a natural map $f\colon \YA^{\Br}\ra \cXA^{\Br}$.
By Lemma \ref{lemm_Br_iso}, $f^*:\Br\cX\to \Br Y$ is an isomorphism.
By the functorial definition of obstructions \eqref{eq_obs},
 we know that for any $x\in \cXA^{\Br}$, any lift $y\in \YA$ along $f$
 must lie in $\YA^{\Br}$.
Since Proposition \ref{prop_Ak_desc} shows that
 the natural map $f\colon \YA\ra \cXA$ is surjective, the induced map
 $\YA^{\Br}\ra \cXA^{\Br}$ is surjective.
}

\lemm{ \label{lemm_pi1_iso}
Let $\cX=[Y/G]$ where $Y$ is a geometrically connected
 $k$-scheme and $G$ is
 a simply connected $k$-group.
Then the natural map  $f: Y\ra \cX$ induces an isomorphism
 $\pi_1^\et(f): \pi_1^\et(Y)\xra{\sim} \pi_1^\et(\cX)$.
}
\pf{
\cite[Thm. A.10]{noohi04fundamental} shows that the quotient map $f: Y \to \cX$ is a fibration (in the sense of
   \cite[Def. A.9]{noohi04fundamental}), there is a homotopy exact sequence of \'{e}tale fundamental groups attached to $f$
   \eqn{
     \pi_1^\et(G_{\bar k}) \to \pi_1^\et(Y)
\xrightarrow{\pi_1^\et(f)}
\pi_1^\et(\cX) \to \pi_0(G_{\bar k}) \to \{*\}.
    }
Since $G$ is simply connected, $G_{\bar k}$ is connected and
simply connected. Thus
\[
\pi_0(G_{\bar k})=1 \text{ and } \pi_1^\et(G_{\bar k})=1 .
\]
Hence $\pi_1^\et(f): \pi_1^\et(Y) \xrightarrow{\sim} \pi_1^\et(\cX)$ is an isomorphism. The proof is complete.
}

\lemm{ \label{lemm_fin_tor_pb}
Let $\cX=[Y/G]$ where $Y$ is a geometrically connected
 $k$-variety and $G$ is
 a simply connected $k$-group.
Let $F$ be a finite $k$-group.
Then the pullback induces a natural equivalence of groupoids
 $\Tors(\cX, F)\xra{\sim} \Tors(Y, F)$, i.e.,
 any $F$-torsor $Z'\xra{F} Y$ is isomorphic to
   the pullback $\cZ\tm_\cX Y$ for some $F$-torsor $\cZ\xra{F} \cX$.
}
\pf{
  For any connected
   algebraic stack $\cY$, let $\Chp_{/\cY}^\finet$ be the ($1$-)category
   of representable finite {\etale} algebraic stacks over $\cY$.
  Let $F_\cY: \Chp_{/\cY}^\finet\ra \Set$
   be the functor sending any $\cZ\ra \cY$
   to the geometric fiber of $\cZ$ over $\cY$. See
   \cite[Def. 4.1]{noohi04fundamental}.
  Then $F_\cY$ induces an equivalence
  \eqn{
    F_\cY:  \Chp_{/\cY}^\finet\xra{\sim} \pi_1^\et(\cY)\text{-}\Set.
  }
  Moreover, for any finite $k$-group $F$,
   under this equivalence, the groupoid
   $\Tors(\cY, F)\subseteq \Chp_{/\cY}^\finet$
   corresponds to the  groupoid whose objects are the
   finite set $F(\ol k)$
   with $\pi_1^\et(\cY)$ acting  via a surjective homomorphism
   $\pi_1^\et(\cY)\ra F(\ol k)$.
  Applying this equivalence to the natural map $f$,  we obtain the
   commutative diagram
  \[
  \begin{tikzcd}[ampersand replacement=\&]
  \Chp_{/\cX}^\finet
    \arrow[r, "\sim" above, "F_\cX" below]
    \arrow[d, "f^*" left]
  \&
  \pi_1^\et(\cX)\text{-}\Set
    \arrow[d, "f^*" right] \\
  \Chp_{/Y}^\finet
    \arrow[r, "\sim" above, "F_Y" below]
  \&
  \pi_1^\et(Y)\text{-}\Set
\end{tikzcd}
  \]
   where the left $f^*$ is the pullback along $f$, and the right $f^*$
   is induced by  $\pi_1^\et(f): \pi_1^\et(Y)\ra \pi_1^\et(\cX)$.
  By Lemma \ref{lemm_pi1_iso},  $\pi_1^\et(f)$  is an isomorphism, it follows
   that both $f^*$ are equivalences.
  In particular, we have an equivalence $f^*: \Tors(\cX, F)\xra{\sim}
   \Tors(\cY, F)$.
 The proof is complete.
}

\lemm{\label{lem_pi0_fin_etale}
Let $V$ be a smooth $k$-variety. Then there exists a finite {\etale}
$k$-scheme $C=\pi_0(V)$ and a $k$-morphism
$q_V:V\ra C$ such that the fibers of
$q_{V,\ol k}:V_{\ol k}\ra C_{\ol k}$ are precisely the connected
components of $V_{\ol k}$. For every $c\in C(k)$, the fiber
$V_c:=V\tm_{C,c} k$ is a smooth geometrically integral $k$-variety.

Moreover, if a geometrically connected $k$-group $H$ acts on $V$, then
$q_V$ is $H$-invariant. Consequently, $q_V$ descends to a morphism
$[V/H]\ra C$, and for every $c\in C(k)$ one has
\eqn{
[V/H]\tm_{C,c} k\simeq [V_c/H].
}
}
\pf{
  By \cite[038E]{sp},
the finite set
$\pi_0(V_{\ol k})$ carries a continuous action of
$\Gal(\ol k/k)$.
By the equivalence between finite {\etale}
$k$-schemes and finite continuous $\Gal(\ol k/k)$-sets
\cite[03QR]{sp}, it corresponds to a finite {\etale}
$k$-scheme $C$.
Over $\ol k$,
send each connected component of
$V_{\ol k}$ to the corresponding point of $C_{\ol k}$. This morphism
is $\Gal(\ol k/k)$-equivariant and therefore descends to a morphism
$q_V:V\ra C$ \cite[040L]{sp}.

For $c\in C(k)$, the geometric fiber $(V_c)_{\ol k}$ is a connected
component of $V_{\ol k}$.
A connected
regular noetherian scheme is integral.
Thus $V_c$ is smooth and
geometrically integral.

Suppose now that $H$ is geometrically connected.
If $D$ is a connected
component of $V_{\ol k}$,
then $H_{\ol k}\tm_{\ol k}D$ is connected,
so its image in the finite discrete set $\pi_0(V_{\ol k})$ consists of
one point.
The identity section of $H_{\ol k}$ shows that this point is
$D$. Hence every geometric connected component is $H_{\ol k}$-stable,
and $q_V$ is $H$-invariant. The assertions for the quotient stack and
its fibers follow by descent and base change.
}

\lemm{\label{lem_rem_complex}
Let $\cV$ be an algebraic stack over $k$, and let $\bfA_k^{\tu{nc}}$
denote the adeles outside the complex places. Define
$\cV(\bfA_k^{\tu{nc}})^{\Br}$ by the same functorial construction as
in \eqref{eq_obs}, with $\bfA_k$ replaced by $\bfA_k^{\tu{nc}}$.
If $\Omega_k^{\CC}$ denotes the set of complex places of $k$, then
\eqn{
\cV(\bfA_k)^{\Br}
=
\cV(\bfA_k^{\tu{nc}})^{\Br}
\tm
\prod_{v\in\Omega_k^{\CC}}\cV(\CC).
}
}
\pf{
This follows directly from the functorial definition \eqref{eq_obs}
and the fact that $\Br(\CC)=0$. Thus the Brauer condition imposes no
condition on the components at complex places.
}

\prop{ \label{prop_etBr_desc}
Let $\cX=[Y/\SL_n]$, where $Y$ is a smooth geometrically integral
 $k$-variety.
Let  $f$ be the natural map  $Y\ra \cX$.
Then for any $x\in \cXA^{\etBr}$, any lift $y\in \YA$ along $f$
 must lie in $\YA^{\etBr}$.
In particular, $f$ induces a surjective map
 $\YA^{\etBr}\ra \cXA^{\etBr}$.
}
\pf{
Let $x\in \cXA^{\etBr}$ and $y\in \YA$ be a lift along $f$.
  We want to show that $y\in \YA^\etBr$.
  By the definition of the $(\etBr)$ obstruction,
   we need to show that
   for any
   finite torsor  $Z'\xra{F} Y$,
   $y$ lifts to $Z'^\s(\bfA_k)^{\Br}$    for some $\s\in H^1(k, F)$,
  where $Z'^\s\xra{F^\s} Y$ is the twist. Lemma
\ref{lemm_fin_tor_pb}  tells us that $Z'$ is isomorphic to
   the pullback $Z=\cZ\tm_\cX Y$ for some finite torsor $\cZ\xra{F} \cX$.
  Since $x\in \cXA^{\etBr}$, by the  definition of  $(\etBr)$ again,
   there is $\s\in H^1(k, F)$ such that
   $x$ lifts along the twist $\cZ^\s\xra{F^\s} \cX$ to
   some $z\in\cZ^\s(\bfA_k)^{\Br}$.
  We know that $Z^\s \xra{\sim} \cZ^\s\tm_\cX Y$
   is the twist of $Z\xra{F} Y$,
   and $Z^\s\ra \cZ^\s$ is an
   $\SL_n$-torsor since $f$ is an $\SL_n$-torsor and
   $H^1(k, \SL_n)$ is trivial (Lemma \ref{lemm_H1_SLn_field}).
  By the universal property of the fiber product, there exists  $z'\in Z^\s(\bfA_k)$
    lifting both $z$ and $y$,  which we depict as follows
  \[
  \begin{tikzcd}[ampersand replacement=\&]
  \Spec \bfA_k
    \arrow[to=2-3, bend left, "z" above]
    \arrow[to=2-2, dotted, "z'"]
    \arrow[to=3-2, bend right, "y" below]
  \& \& \\
  \& Z^\s
    \arrow[r, "\SL_n" above]
    \arrow[d, "F^\s" left]
  \& \cZ^\s
    \arrow[d, "F^\s" right] \\
  \& Y
    \arrow[r, "\SL_n" above, "f" below]
  \& \cX
\end{tikzcd}
  \]
Since $Y$ is a smooth $k$-variety,
 so is $Z^\s$.
Let
$\pi: Z^\s\ra C:=\pi_0(Z^\s)$ be the geometric connected-component
morphism constructed in Lemma~\ref{lem_pi0_fin_etale}.
Since $\SL_n$ is geometrically connected,
$\pi$ is $\SL_n$-invariant.
As
$\cZ^\s\simeq[Z^\s/\SL_n]$,
the morphism $q$ descends to
\eqn{
\ol \pi:\cZ^\s\ra C.
}

Let $z^{\tu{nc}}$ denote the components of $z$ outside the complex
places.
By Lemma~\ref{lem_rem_complex},
$z^{\tu{nc}}\in\cZ^\s(\bfA_k^{\tu{nc}})^{\Br}$;
hence, by
functoriality,
\eqn{
\ol \pi (z^{\tu{nc}})\in C(\bfA_k^{\tu{nc}})^{\Br}.
}
Since $C$ is finite {\etale} over $k$, \cite[Prop. 3.3]{lx15very}
gives
\eqn{
C(\bfA_k^{\tu{nc}})^{\Br}=C(k).
}
Thus there exists $c\in C(k)$
such that the diagram
\eqn{ \xymatrix{
  \Spec(\bfA_k^{\tu{nc}}) \ar[r]^-q \ar[d]^-{z^{\tu{nc}}}
     &\Spec k \ar[d]^-{c} \\
   \cZ^\s \ar[r]^-{\ol \pi} & C
}}
 commutes.

Set
\eqn{
Z_c^\s:=Z^\s\tm_{C,c}k,
\qquad
\cZ_c^\s:=\cZ^\s\tm_{C,c}k.
}
By Lemma \ref{lem_pi0_fin_etale},  $Z_c^\s$ is a smooth
geometrically integral $k$-variety and
$\cZ_c^\s\simeq[Z_c^\s/\SL_n]$. Moreover, $z^{\tu{nc}}$ factors
through $\cZ_c^\s$, and, since $z'$ maps to $z$,
$z'^{\tu{nc}}$ factors through $Z_c^\s$.

The section $c:\Spec k\ra C$ is open and closed, so
$\cZ_c^\s$ is an open-and-closed substack of $\cZ^\s$. Hence the
restriction map $\Br(\cZ^\s)\ra\Br(\cZ_c^\s)$ is surjective. By the
functorial definition \eqref{eq_obs}, we therefore have
\eqn{
z^{\tu{nc}}\in\cZ_c^\s(\bfA_k^{\tu{nc}})^{\Br}.
}
For each complex place $v$, choose a point
$\widetilde z'_v\in Z_c^\s(\CC)$ and let $\widetilde z_v$ be its
image in $\cZ_c^\s(\CC)$. Together with the original non-complex
components $z^{\tu{nc}}$, these choices give adelic points
$\widetilde z'\in Z_c^\s(\bfA_k)$ and
$\widetilde z\in\cZ_c^\s(\bfA_k)$, with $\widetilde z'$ lifting
$\widetilde z$. Lemma~\ref{lem_rem_complex} gives
$\widetilde z\in\cZ_c^\s(\bfA_k)^{\Br}$.

We can now apply Proposition~\ref{prop_Br_desc} to the $\SL_n$-torsor
$Z_c^\s\ra\cZ_c^\s$, because $Z_c^\s$ is smooth and geometrically
integral. It follows that
$\widetilde z'\in Z_c^\s(\bfA_k)^{\Br}$. Applying
Lemma~\ref{lem_rem_complex} once more yields
\eqn{
z'^{\tu{nc}}=\widetilde z'^{\,\tu{nc}}
\in Z_c^\s(\bfA_k^{\tu{nc}})^{\Br}.
}
By functoriality for $Z_c^\s\ra Z^\s$, we obtain
$z'^{\tu{nc}}\in Z^\s(\bfA_k^{\tu{nc}})^{\Br}$.
Applying Lemma~\ref{lem_rem_complex} to the original complex
components of $z'$ yields $z'\in Z^\s(\bfA_k)^{\Br}$.

Thus $y$ lifts to $Z^\s(\bfA_k)^{\Br}\simeq Z'^\s(\bfA_k)^{\Br}$
for a suitable $\s\in H^1(k,F)$. Since the finite torsor $Z'\ra Y$
was arbitrary, the first statement
   of the proposition has been shown.

   % FIXME: \cred{it is geo. int. ? SOLVED
Since any $x\in\cXA$ can lift to some $y\in \YA$ by Proposition \ref{prop_Ak_desc},
   combining the first statement, we know that $f$ induces a surjective map
   $\YA^{\etBr}\ra \cXA^{\etBr}$.
  The proof is complete.
}

\section{Comparing relations} \label{cmp}

In this section, we are going to show the following relations of
 various obstructions on quotient stacks $[X/G]$.
In particular,  $(\detBr)$ is the finest among
 $\desc$, $\conn$, $(\ddesc)$, $(\etBr)$, $\sdesc$, $(\fdesc)$
  and $\Br$.
 % for iterated precompositions by $\desc$,
 %   FIXME: geo int can NOT passing to a G-torsor !!!

\thm{ \label{thm_rel}
Let $\cX$ be an algebraic $k$-stack.
\enmt{[\upshape (i)]
\item Then we have the following relations
\[
\begin{tikzcd}[ampersand replacement=\&]
  \cXA^\detBr
    \arrow[r, equal, "\arrowlabel{deB=eB}" above] \&
  \cXA^\etBr
    \arrow[r, hookrightarrow, "\arrowlabel{eB<=dd}" above]
    \arrow[d, hookrightarrow, "\arrowlabel{eB<=B}" left] \&
  \cXA^\ddesc
    \arrow[r, hookrightarrow, "\arrowlabel{dd<=fd}" above] \&
  \cXA^\fdesc
    \arrow[r, equal, "\arrowlabel{fd=d}" above] \&
  \cXA^\desc
    \arrow[d, hookrightarrow, "\arrowlabel{d<=c}" left] \\
  \&
  \cXA^{\Br}
    \arrow[r, equal, "\arrowlabel{B=sd}" above] \&
  \cXA^\sdesc
    \arrow[rr, hookrightarrow, "\arrowlabel{sd<=c}" above] \& \&
  \cXA^\conn
\end{tikzcd}
\]
 where $\cXA^\detBr$ is defined to be
\eqn{
\cXA^\detBr=\bigcap_{f:Y\xra{G} \cX \text{ torsor under
 linear $k$-group $G$}} \bigcup_{\s\in H^1(k, G)}
 f^\s(Y^\s(\bfA_k)^{\etBr}).
}
 For \arrowref{B=sd} and \arrowref{sd<=c}, assume that $\cX$ is smooth
  of finite type and is either {\DM} or
 Zariski-locally the quotient of a smooth geometrically integral
 $k$-variety by a linear $k$-group. For \arrowref{deB=eB} and \arrowref{eB<=dd}, assume that $\cX=[X/G]$, where $X$ is a
   smooth quasi-projective geometrically integral
 $k$-variety and
 $G$ is a linear $k$-group.

\item \label{it_all_func}
All obstructions appearing in the diagram are functorial.
}}
\pf{
The functoriality \eqref{it_all_func}
  follows from Remark \ref{rk_functorial}
 and Lemma \ref{lemm_func} \eqref{it_func}.

The inclusions \arrowref{dd<=fd} and \arrowref{d<=c}
 are trivial.
The inclusion \arrowref{eB<=B} is due to Lemma
 \ref{lemm_func} \eqref{it_cap}.
The equality \arrowref{fd=d} is \cite[Cor. 4.5]{wl24stackdd} and
 \arrowref{B=sd} is \cite[Thm. 5.2]{lh23stackbm}.
The inclusion \arrowref{sd<=c} is \cite[Thm. 5.4]{lh23stackbm}.
It suffices to show \arrowref{deB=eB} and \arrowref{eB<=dd}, which is
 done by the following two lemmas.
}

\lemm{ \label{lemm_detBr=etBr}
Let $\cX=[X/G]$ where $X$ is a smooth quasi-projective
  geometrically integral
 $k$-variety and
 $G$ is a linear $k$-group.
Then $\cXA^{\detBr} = \cXA^{\etBr}$.
}
\pf{
  The inclusion $\cXA^{\detBr} \subseteq \cXA^{\etBr}$ easily follows from Lemma
   \ref{lemm_func} \eqref{it_cap},
   and then we want to show the converse.
  By Lemma \ref{lemm_quo_SLn},
   we may assume that   $\cX=[Y/\SL_n]$ where
   $Y = \SL_n\tm_k^G X$ = $(\SL_n\tm_k X) /G$ is a smooth quasi-projective geometrically
    integral $k$-variety.
  By  \cite[Thm. 1.1]{cao20sous}  we have, for any linear $k$-group $H$
   and any torsor $f: Z\xra{H} Y$,
  \eqn{
    \YA^\etBr = \bigcup_{\s\in H^1(k, H)}f^\s( Z^\s(\bfA_k)^\etBr).
  }
  Intersecting over all such linear groups and torsors, we obtain
   that $\YA^\etBr=\YA^\detBr$.
  Next, by Proposition \ref{prop_etBr_desc},  the natural map $p: Y\ra \cX$ induces a surjection $\YA^{\etBr}\ra \cXA^{\etBr}$.
  It follows that
  \eqn{
    \cXA^\etBr = p(\YA^\etBr) = p(\YA^\detBr)\subseteq \cXA^\detBr,
  }
   where the inclusion follows from functoriality
    (Theorem \ref{thm_rel} \eqref{it_all_func}).
  The proof is complete.
}

\lemm{ \label{lemm_etBr<=ddesc}
Let $\cX=[X/G]$ where $X$ is a smooth
 quasi-projective geometrically integral
 $k$-variety and
 $G$ is a linear $k$-group.
Then $\cXA^{\etBr}\subseteq \cXA^{\ddesc}$.
}
\pf{
  By Lemma \ref{lemm_quo_SLn},
   we may assume that   $\cX=[Y/\SL_n]$ where
   $Y = \SL_n\tm_k^G X$ = $(\SL_n\tm_k X) /G$ is a smooth quasi-projective geometrically
    integral $k$-variety.
  By  \cite[Thm. 1.2]{cao20sous} and \cite[Thm. 1.5]{cdx19comparing} we have
  \eqn{
    \YA^\etBr = \YA^\desc = \YA^\ddesc.
  }
  Next, by Proposition \ref{prop_etBr_desc},  the  map
   $\YA^{\etBr}\ra \cXA^{\etBr}$ is surjective,
   which is induced  by the natural map $f: Y\ra \cX$.
  It follows that
  \eqn{
    \cXA^\etBr = f(\YA^\etBr) = f(\YA^\ddesc)\subseteq \cXA^\ddesc,
  }
   with the inclusion by functoriality
    (Theorem \ref{thm_rel} \eqref{it_all_func}).
  The proof is complete.
}

%bib
\bibliography{./unibib}
\bibliographystyle{amsalpha}
\end{document}